\documentclass[3p,preprint]{elsarticle}

\usepackage{amsmath}
\usepackage{amsfonts}
\usepackage{amssymb,latexsym}
\usepackage{mathrsfs}
\usepackage{amsthm}
\usepackage{graphicx}
\usepackage {xcolor}

\newtheorem{Corollary}{Corollary}

\newtheorem{Remark}{Remark}
\newtheorem{theorem}{Theorem}

\usepackage{color}
\usepackage{newlfont}
\usepackage{subfigure}
\usepackage{verbatim}
\usepackage{rotating}
\usepackage{multirow}
\usepackage{units}
\usepackage{bm}
\usepackage[english]{babel}
\usepackage[utf8]{inputenc}
\usepackage{algorithm}
\usepackage[noend]{algpseudocode}
\usepackage{booktabs}
\usepackage{caption}
\usepackage{hyperref}

\journal{}

\begin{document}

\begin{frontmatter}

\title{Compressed QMC volume and surface integration on union of balls}



\author[address-PD]{Giacomo Elefante}
\ead{giacomo.elefante@unipd.it}

\author[address-PD]{Alvise Sommariva}{\corref{mycorrespondingauthor}}
\cortext[mycorrespondingauthor]{Corresponding author: Alvise Sommariva}
\ead{alvise@math.unipd.it}

\author[address-PD]{Marco Vianello}
\ead{marcov@math.unipd.it}

\address[address-PD]{University of Padova, Italy}


\begin{abstract}
We discuss an algorithm for Tchakaloff-like compression of Quasi-MonteCarlo (QMC) volume/surface integration on union of balls (multibubbles). The key tools are   
Davis-Wilhelmsen theorem on the so-called “Tchakaloff sets” for positive linear functionals on polynomial spaces, and Lawson-Hanson algorithm for NNLS. We provide the corresponding Matlab package together with several examples. 
\end{abstract}



\end{frontmatter}

\section{Introduction}
Numerical modelling by finite collections of disks, balls and spheres is relevant within different application fields. Problems involving intersection, union and difference of such geometrical objects arise for example in molecular modelling, computational geometry, computational optics, wireless network analysis; cf., e.g., \cite{ABI88,BX10,DQS20,KN22,LVZ14,P13} with the references therein. 
A basic problem is the computation of areas and volumes of such sets, followed by the 
more difficult task of computing volume and surface integrals there by suitable quadrature formulas.

Indeed, the numerical quadrature problem on intersection and union of planar disks has been recently treated in \cite{SV17,SV22}, providing low-cardinality algebraic formulas with positive weights and interior nodes.
On the other hand, though there is some literature, mainly in the molecular modelling field, 
on the computation of volumes and surface areas 
of arbitrary union of balls (multibubbles), to our knowledge specific numerical integration codes on such domains are not available yet.  

\begin{figure}
    \centering
    \subfigure[]
        {\includegraphics[height=1.81in]{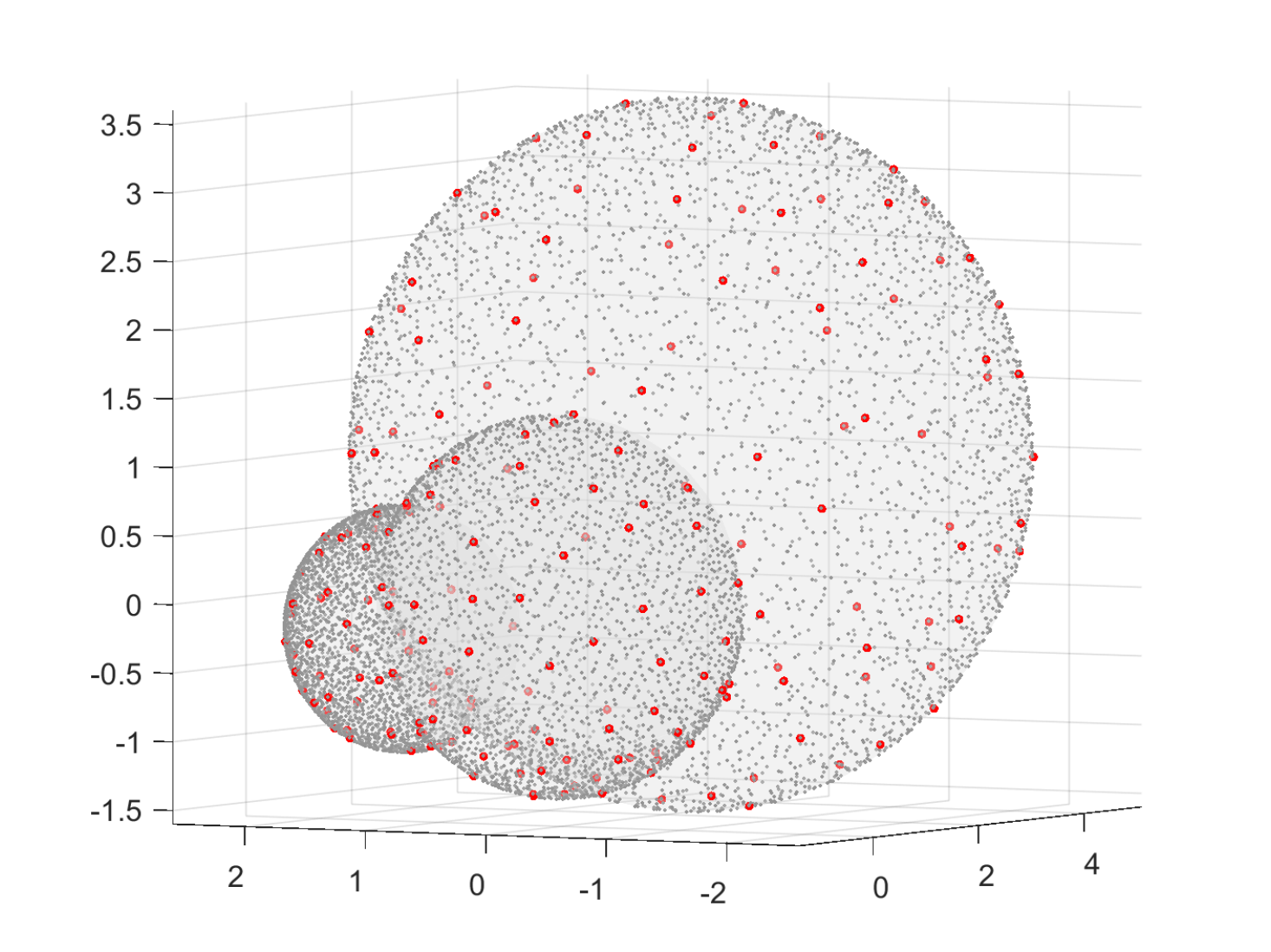}}
     ~ 
    \subfigure[]
        {\includegraphics[height=1.81in]{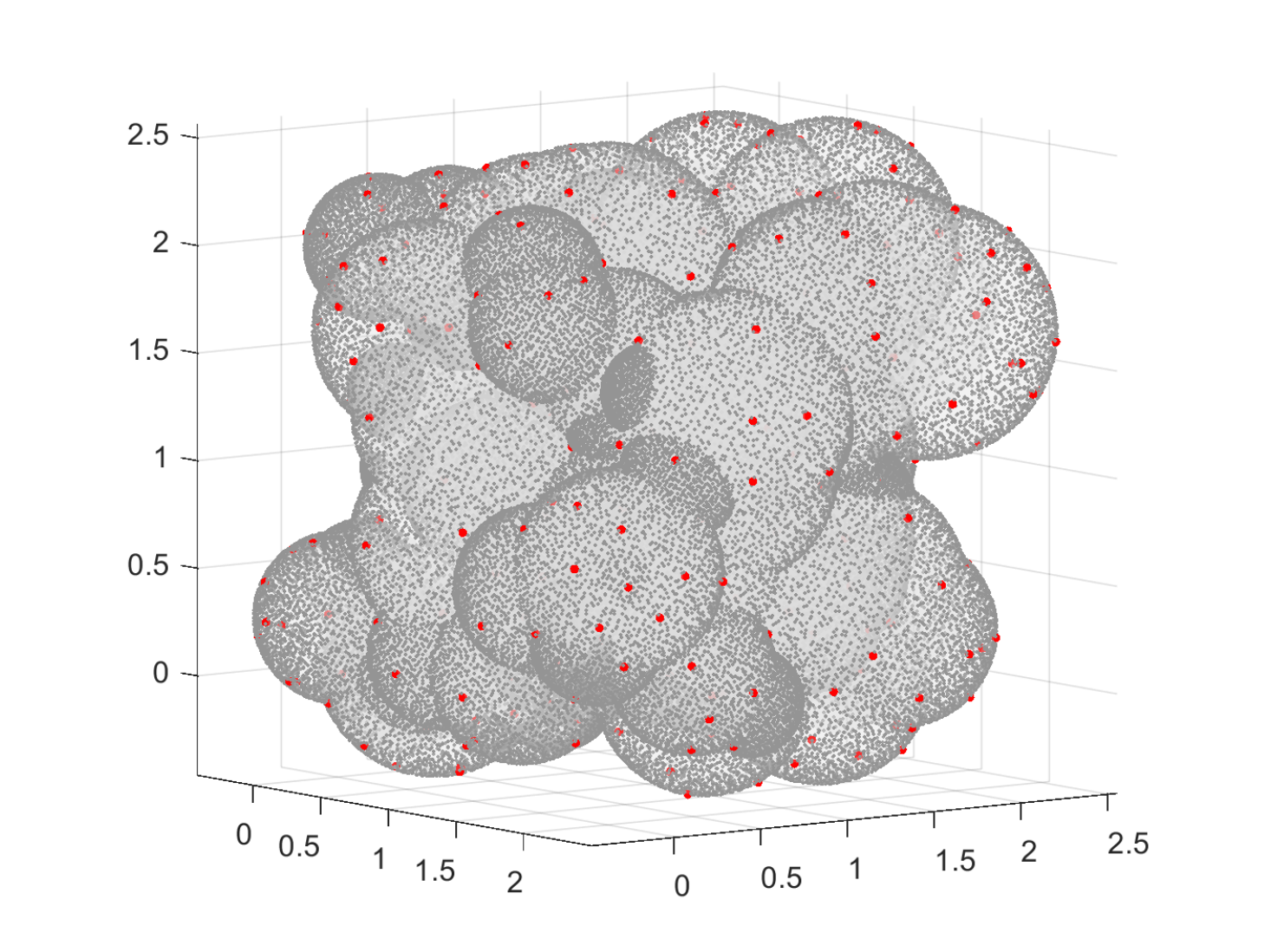}}
    \caption{Compressed QMC points (red) extracted from low-discrepancy points (grey) on the surface of ball union at degree $n=9$. Left: 200 points extracted from about 8200 (3 balls), compression ratio 43; Right: 220 points extracted from about 69000 (100 balls), compression ratio over 300.}
    \label{compressed-deg9}
\end{figure}

In this paper, we begin to fill the gap by providing compressed Quasi-Montecarlo (QMC) formulas 
for volume and surface integration on multibubbles, along the lines of \cite{ESV22}. Such formulas {\em preserve the approximation power of QMC} up to the best uniform polynomial approximation error 
of a given degree to the integrand, but using a {\em much lower number of sampling points}; see Figure \ref{compressed-deg9} for two examples of multibubbles and QMC sampling compression.
The key tools are Davis-Wilhelmsen theorem on the so-called ``Tchakaloff sets'' for positive linear functionals and Lawson-Hanson algorithm for 
NNLS, which allows to extract a set of ``equivalent'' re-weighted nodes from a huge 
low-discrepancy sequence. 

{ We stress that differently from \cite{ESV22}, the present approach is able to compress not only QMC volume integration, but also QMC integration on compact subsets of algebraic surfaces (in particolar, the surface of a multibubble which is a subset of a union of spheres). Notice that one of the main difficulties in surface instances, consists in adapting the compression algorithm
to work on spaces of polynomials restricted to an algebraic variety, finding an appropriate polynomial basis. Indeed, to our knowledge the present work is the first attempt in this direction within the QMC framework.} 

The paper is organized as follows. 
In Section 2, we discuss theoretical and computational issues of QMC compression for volume and surface integration in $\mathbb{R}^3$. In Section 3 we describe our implementation on 3d multibubbles, presenting several numerical tests. The open-source codes 
are freely available at \cite{ESV23}.

\section{Compressed QMC formulas}

Compression of QMC formulas is
nothing but a special instance of {\em discrete measure compression}, a topic which has received an increasing attention in the literature of the last decade, in both the probabilistic and the deterministic setting. 
Indeed, several papers and some software 
have been devoted to the extraction of a smaller set of re-weighted mass points from the support of a high-cardinality discrete measure, with the constraint 
of preserving its moments up to a given 
polynomial degree; cf., e.g., \cite{H21,LL12,PSV17,SV15,Tche15} with the references therein. 

From the quadrature point of view, this topic has a strong connection with the famous Tchakaloff 
theorem \cite{Tcha57} on the existence of low-cardinality formulas with positive weights. 
On the other hand, Tchakaloff theorem itself 
is contained in a somewhat deeper but somehow overlooked result by Wilhelmsen \cite{W76} on the the discrete representation of positive linear functionals on finite-dimensional 
function spaces (which generalizes a previous result by 
Davis \cite{D67}). Indeed, only quite recently this theorem has been rediscovered as a basic tool for positive cubature via adaptive NNLS moment-matching, cf. \cite{ESV22,L21,SV21,SV23}. 

\begin{theorem} (Davis, 1967 - Wilhelmsen, 1976) Let $\Psi$ be the linear span of continuous, 
real-valued, linearly independent functions $\{\phi_j\}_{j=1,\ldots,N}$ 
defined on a compact set $\Omega \subset {\mathbb{R}}^d$. 
Assume that $\Psi$ satisfies the Krein condition (i.e. there is at least 
one $f \in \Psi$ which does not vanish on $\Omega$) and that $L$ is a positive 
linear functional on $\Psi$, i.e. $L(f)>0$ for every $f\in \Psi$, $f \geq 0$ not vanishing everywhere in $\Omega$. 

If $\{P_i\}_{i=1}^{\infty}$ is an everywhere dense subset of $\Omega$, 
then for sufficiently large $m$, the set 
$X_m=\{P_i\}_{i=1,\ldots,m}$ is 
a Tchakaloff set, i.e. there exist weights $w_k > 0$, $k=1,\ldots,\nu$, and nodes $\{\mathcal{T}_k\}_{k=1,\dots,\nu}
\subset X_m \subset \Omega$, 
with $\nu={\mbox{card}}(\{\mathcal{T}_k\}) \leq N$, such that 
\begin{equation} \label{tchq}
L(f)=\sum_{k=1}^\nu w_k f(\mathcal{T}_k)\;,\;\;\forall f \in \Psi\;.
\end{equation}
\end{theorem}

\vskip0.5cm

As an immediate consequence, we may state the following

\begin{Corollary}
Let $\lambda$ be a positive measure on $\Omega$, such that $\mbox{\em{supp}}(\lambda)$ is determining for 
$\mathbb{P}_n^d(\Omega)$, the space of total-degree polynomials of degree not exceeding $n$, restricted to $\Omega$ (i.e., a polynomial in $\mathbb{P}_n^d(\Omega)$ vanishing there 
vanishes everywhere on $\Omega$). Then the thesis of Theorem 1 holds for $L(f)=\int_\Omega{f\,d\lambda}$.
\end{Corollary}

Indeed, the integral of a nonnegative and not everywhere vanishing polynomial $f\in \mathbb{P}_n^d(\Omega)$ must be positive (otherwise $f$ would vanish on ${\mbox{\em supp}}(\lambda)$).
Observe that the classical version of Tchakaloff theorem corresponds to 
$$
L(f)=L_{\mbox{\tiny{INT}}}(f)=\int_\Omega{f(P)\,dP}\;,
$$
with $\Psi=\mathbb{P}_n^d(\Omega)$ and 
\begin{equation} \label{dim}
N=N_n^d=dim(\mathbb{P}_n^d(\Omega))\;.
\end{equation}

From now on we shall concentrate on the 3-dimensional case ($d=3$), though most considerations could be extended in general dimension. 
Notice that 
the formulation of Davis-Wilhelmsen theorem is sufficiently general to include volume integrals,  i.e. $\Omega$ is the closure of a bounded open set and $N=dim(\mathbb{P}_n^d(\mathbb{R}^3))={n+3 \choose 3}=(n+1)(n+2)(n+3)/6$, as well as surface integrals on compact subsets of an algebraic variety (in this case $dP=d\sigma$ for the surface measure). In the latter case the dimension of the polynomial space could collapse,  
for example with $\Omega=S^2\subset \mathbb{R}^3$ we have $N=(n+1)^2<{n+3 \choose 3}=(n+1)(n+2)(n+3)/6$.

On the other hand, Wilhelmsen theorem can also 
be applied to a discrete functional like a QMC formula applied to $f\in C(\Omega)$ 
\begin{equation} \label{QMC}
L(f)=L_{\mbox{\tiny{QMC}}}(f)=\frac{\mu(\Omega)}{M}\,\sum_{i=1}^M{f(P_i)}\approx 
\int_\Omega{f(P)\,dP}\;,
\end{equation}
where 
$$X_M=\{P_i\}_{i=1,\ldots,M}\;,\;\;M>N\;,$$ 
is a {\em low-discrepancy sequence} on $\Omega$, and $\mu(\Omega)$ can be  either a volume or a surface area. 
Typically one generates a low-discrepancy sequence of cardinality say $M_0$ on a bounding box or bounding surface $\mathcal{B}\supseteq \Omega$, from which the low-discrepancy sequence on $\Omega$ is extracted by a suitable in-domain algorithm.  We observe that if $\mu(\Omega)$ is unknown or difficult to compute (as in the case of multibubbles), it can be approximated  as $\mu(\Omega)\approx \mu(\mathcal{B})M/M_0$.

Positivity of the functional for $f\in \Psi=\mathbb{P}_n^3(\Omega)$ is ensured whenever the set $X_M$ is $\mathbb{P}_n^3(\Omega)$-determining, i.e. polynomial vanishing there vanishes everywhere on $\Omega$, 
or equivalently $dim(\mathbb{P}_n^3(X_M))=N=dim(\mathbb{P}_n^3(\Omega))$, or even 
\begin{equation} \label{ranks}
rank(V_M)=N\;,
\end{equation}
where 
\begin{equation} \label{vdm}
V_M=V^{(n)}(X_M)=[\phi_j(P_i)]\in \mathbb{R}^{M\times N}
\end{equation}
is the corresponding rectangular Vandermonde-like matrix. Notice that, $X_M$ being a sequence, for every $k\leq M$ we have that 
\begin{equation} \label{submatrices}
V_k=V^{(n)}(X_k)=[(V_M)_{ij}]\;,\;\;1\leq i \leq k\;,\;1\leq j\leq N\;.
\end{equation}
The full rank requirement for $V_M$ is not restrictive, in practice. 
Indeed, the probability that $det(V_N)=0$ dealing with uniformly distributed points is null, since the former equation defines the zero set of a polynomial in $\Omega^N$, whose product measure is null (cf., e.g.,   \cite[\S\S 3-4]{PPL21} for a more complete discussion on this point).  

By Theorem 1, when $M\gg N$ we can then try to find a Tchakaloff set $X_m$, with $N\leq m<M$, such that a sparse nonnegative solution vector $u$ exists to the underdetermined {\em moment-matching system}
\begin{equation} \label{mom-match}
V_m^t u=p=V_M^t e\;,\;\;e=\frac{\mu(\Omega)}{M}\,(1,\dots,1)^t\;.
\end{equation}

In practice, we solve (\ref{mom-match}) via Lawson-Hanson active-set method \cite{LH95} applied to the NNLS problem 
\begin{equation} \label{NNLS}
\min_{u\geq 0} \|V_m^t u-p\|_2\;, 
\end{equation}
accepting the solution when the residual size is small, say 
\begin{equation} \label{res-tol}
\|V_m^t u-p\|_2<\varepsilon
\end{equation}
where $\varepsilon$ is a given tolerance. The nonzero components of $u$ then determine the nodes and weights of a compressed QMC formula {extracted from $X_m$}, that is $\{w_k\}=\{u_i:\,u_i>0\}$ and $\{\mathcal{T}_k\}=\{P_i:\,u_i>0\}$, giving 
\begin{equation} \label{compressed}
L^\ast_{\mbox{\tiny{QMC}}}(f)=\sum_{k=1}^\nu w_k f(\mathcal{T}_k)\;,\;\;\nu \leq N\ll M\;\;,
\end{equation} 
where $L^\ast_{\mbox{\tiny{QMC}}}(f)=L_{\mbox{\tiny{QMC}}}(f)$ for every $f\in \mathbb{P}_n^3(\Omega)$.

Notice that existence of a representation like (\ref{compressed}) for $m=M$ is ensured by Caratheodory theorem on finite-dimensional conic combinations,  
applied to the columns of $V_M^t$ (cf. \cite{PSV17} for a full discussion on this point in the general framework of discrete measure compression). In such a way, however, we would have to work with a much larger matrix, that is we would have to solve directly 
\begin{equation} \label{big-NNLS}
\min_{u\geq 0} \|V_M^t u-p\|_2\;. 
\end{equation}
On the contrary, solving (\ref{NNLS}) on an increasing sequence of smaller problems $m:=m_1,m_2,m_3,\dots$ with $m_1<m_2<m_3<\dots\leq M$,
\begin{equation} \label{small-NNLS}
\min_{u\geq 0} \|V_{m_j}^t u-p\|_2\;,\;\;j=1,2,3,\dots\;,\;m_1\geq N\;,
\end{equation}
corresponding to increasingly dense subsets $X_{m_1}\subset X_{m_2} \subset \dots \subseteq  X_M$
(say, ``bottom-up''), until the residual becomes sufficiently small, could substantially lower the computational cost. Indeed, as shown in \cite{ESV22}, 
with a suitable choice of the sequence $\{m_j\}$ the residual becomes extremely small in few iterations, with a final extraction cardinality much lower than $M$.

Concerning the approximation power of QMC compression, following \cite{ESV22} it is easy to derive the following error estimate 
$$
|L^\ast_{\mbox{\tiny{QMC}}}(f)-L_{\mbox{\tiny{INT}}}(f)|
\leq \mathcal{E}_{\mbox{\tiny{QMC}}}(f) 
+2\,\mu(\Omega)\,E_n(f;X)
$$
\begin{equation} \label{equiv}
\leq \mathcal{E}_{\mbox{\tiny{QMC}}}(f) 
+2\,\mu(\Omega)\,E_n(f;\Omega)\;,
\end{equation}
valid for every $f\in C(\Omega)$, where 
$\mathcal{E}_{\mbox{\tiny{QMC}}}(f)=|L_{\mbox{\tiny{QMC}}}(f)-L_{\mbox{\tiny{INT}}}(f)|$ and 
we define $E_n(f;K)=\inf_{\phi\in \mathbb{P}_n^3(K)}
{\|f-\phi\|_{\infty,K}}$ 
with $K$ discrete or continuous compact set.

The meaning of (\ref{equiv}) is that the compressed QMC functional $L^\ast_{\mbox{\tiny{QMC}}}$ retains the approximation power of the original QMC formula, up to a quantity proportional to the best polynomial approximation error to $f$ in the uniform norm on $X$ (and hence by inclusion in the uniform norm on $\Omega$). We recall that the latter can be estimated depending on the regularity of $f$ by multivariate Jackson-like theorems, cf. e.g. \cite{P09} for volume integrals where $\Omega$ is the closure of a bounded open set. 

On the other hand, we do not deepen here 
the vast and well-studied topic of QMC convergence and error estimates, recalling only that (roughly) the QMC error $\mathcal{E}_{\mbox{\tiny{QMC}}}(f)$ is close to $\mathcal{O}(1/M)$ for smooth functions, to be compared with the $\mathcal{O}(1/\sqrt{M})$ error of MC. For basic concepts and results of QMC theory like   
discrepancy, star-discrepancy, 
Hardy-Krause variation, Erd\"{o}s-Tur\'{a}n-Koksma and Koksma-Hlawka inequalities, we refer the reader to devoted surveys like e.g. \cite{DP10}.

\begin{Remark} \label{QMC-union-seq}
{\em 
The QMC compression algorithm can be easily extended to the case where $\Omega$ (either a volume or a surface) is the finite union of nonoverlapping subsets, say $\Omega=\cup_{\ell=1}^L\Omega_\ell$, such that sequences of low-discrepancy points are known on bounding sets $\mathcal{B}_\ell\supset \Omega_\ell$. In this case the overall QMC points are $X=\cup_{\ell=1}^L Y_\ell$, with $Y_\ell=\{P_{\ell, i}\}_{i=1}^{M_\ell}$ and $M=card(X)=\sum_{\ell=1}^L{M_\ell}$, 
where $Y_\ell$ are the low-discrepancy points of $\mathcal{B}_\ell$ lying in $\Omega_\ell$. We stress that the low-discrepancy points have to be chosen alternatively in order to construct an evenly distributed sequence $X_M$ on the whole $\Omega$, picking the first point in each $\Omega_\ell$, then the second point in each $\Omega_\ell$ and so on, i.e. the sequence 
$\{P_{1,1},P_{2,1},\ldots,P_{L,1},P_{1,2},P_{2,2},\ldots,P_{L,2},\ldots\}$.

Moreover, by additivity of the integral the QMC functional 
becomes 
\begin{equation} \label{QMC-union}
L_{QMC}(f)=\sum_{\ell=1}^L{\sum_{i=1}^{M_\ell}w_{\ell,i}f(P_{\ell,i})}\approx
\sum_{\ell=1}^L{\int_{\Omega_\ell}{f(P)\,dP}}=\int_{\Omega}{f(P)\,dP}\;,
\end{equation}
where $w_{\ell i}=\mu(\Omega_\ell)/M_\ell$, $i=1,\dots,M_\ell$, and hence the QMC moments in (\ref{mom-match}) have to be computed with such weights.
}
\end{Remark}

\section{Implementation and numerical tests}

In order to show the effectiveness of the bottom-up compression procedure described in the previous section, we briefly sketch a possible implementation and we present some numerical tests for both, volume and surface integration on arbitrary union of balls (multibubbles). 

Indeed, we compare ``Caratheodory-Tchakaloff'' compression of multivariate discrete measures as implemented 
in the general-purpose package {\em dCATCH} \cite{DMV20-2}, with the bottom-up approach. 
All the tests have been performed with a CPU AMD Ryzen 5 3600 with 48 GB of RAM, running Matlab R2022a.
The Matlab codes and demos, collected in a package named {\em Qbubble}, are freely available at 
\cite{ESV23}.

Below, we first give some highlights on the main features of the implemented algorithm on multibubbles. These are essentially: 
\begin{itemize}
\item for multibubble volume integrals we simply take Halton points of the smaller bounding box $$[a_1,b_1]\times [a_2,b_2]\times [a_3,b_3]\supset \Omega $$ and select those belonging to $\Omega$; for multibubble surface integrals we follow the procedure sketched in Remark {\ref{QMC-union-seq}}, taking on each sphere $\mathcal{B}_\ell$ low-discrepancy mapped Halton points by an area preserving transformation (see (\ref{mapped}) in Section 3.2 below), and then selecting those belonging to the surface;

\item in view of extreme ill-conditioning of the standard monomial basis, we start from the product Chebyshev total-degree basis of the smaller bounding box for $\Omega$ (for either volumes or surfaces), namely 
$$
p_j(x,y,z)=T_{\alpha_1(j)}\left(\sigma_1(x)\right)
T_{\alpha_2(j)}\left(\sigma_2(y)\right)
T_{\alpha_3(j)}\left(\sigma_3(z)\right)\;,\;\;j=1,\dots,J\;,
$$
where 
$J=(n+1)(n+2)(n+3)/6$, $\sigma_i(t)=\frac{2t-b_i-a_i}{b_i-a_i}$, $i=1,2,3$, and 
$j\mapsto \alpha(j)$ corresponds to the graded lexicographical ordering of the triples $\alpha=(\alpha_1,\alpha_2,\alpha_3)$, $0\leq \alpha_1+\alpha_2+\alpha_3\leq n$;

\item for surface integrals we determine a suitable polynomial basis by computing 
the rank and then possibly performing a column selection by QR factorization with column pivoting of the trivariate Chebyshev-Vandermonde matrix; 

\item in order to cope ill-conditioning of the Vandermonde-like matrices $V_{m_j}$ (that increases with the degree), we perform a single QR factorization with column pivoting 
$V_{m_j}=Q_{m_j}R_{m_j}$ to construct 
an orthogonal polynomial basis w.r.t. the discrete scalar product $\langle f,g \rangle_{X_{m_j}}=\sum_{i=1}^{m_j} f(P_i) g(P_i)$ and substitute $V_{m_j}$ by $Q_{m_j}$ in (\ref{small-NNLS}); consequently the QMC moments $p$ in (\ref{mom-match}) have to be modified into $(R^{-1}_{m_j})^tp$ (via Gaussian elimination);

\item the (modified) bottom-up NNLS problems (\ref{small-NNLS}) are solved by the recent 
 implementation of Lawson-Hanson active-set method named LHDM, based on the concept of ``Deviation Maximization'' instead of ``column pivoting'' for the underlying QR factorizations, since it gives experimentally a speed-up of at least 2 with respect to the standard Matlab function {\tt lsqnonneg} (cf. \cite{DODM22,DM22,DMV20}).

\end{itemize}

In the next subsections we present several numerical tests, to show the effectiveness of the bottom-up approach for volume and surface QMC compression on multibubbles.

 \subsection{Volume integration on multibubbles}   

In this subsection we consider volume integration on union of balls (solid multibubbles), namely
\begin{equation} \label{solid}
\Omega=\bigcup_{j=1}^s{B(C_j,r_j)}
\end{equation}
where $B(C_j,r_j)\subset \mathbb{R}^3$ is the closed 3-dimensional ball with center $C_j$ and radius $r_j$.
Here we generate a sequence of Halton points in the smallest Cartesian bounding 
box for $\Omega$ and, then, we select those belonging 
to the union, say $X=\{P_i\}$, simply by checking that $\|P_i-C_j\|_2\leq r_j$ for some $j$. 

More precisely, we consider the following 
\begin{itemize}
\item 
first example: union of the 3 balls with centers $C_1=(0,0,0)$, $C_2=(0,1.3,-0.2)$, $C_3=(2.5,0,1)$ and radii 
$r_1=1.4$, $r_2=0.9$, $r_3=1$, respectively; 
\item 
second example: union of 100 balls 
with randomly chosen and then fixed centers in $[0,2]^3$ and radii in $[0,2,0.6]$. 
\end{itemize}

The results concerning application of the bottom-up approach are collected in Tables {\ref{3balls}}-{\ref{100balls}}, where we compress 
QMC volume integration by more than one million of Halton points, preserving polynomial moments up to 
degree $3,6,9,12,15$ (the moments correspond to the product Chebyshev basis of the minimal Cartesian bounding box for the ball union). 

We start from 2,400,000 Halton points in the bounding box and we set $m_1=2N$ 
and $m_{j+1}=2m_j$, $j\geq 1$. The residual tolerance is $\varepsilon=10^{-10}$.
The comparisons of the present bottom-up compression algorithm, for short 
$Q_c^{bu}$, are made with a global compression algorithm 
that works on the full Halton sequence $X_M$, namely the general purpose discrete measure compressor $dCATCH$ developed in \cite{DMV20-2}, 
which essentially solves directly (\ref{big-NNLS}) by Caratheodory-Tchakaloff subsampling as proposed in \cite{SV15,PSV17}. 

In particular, we display the cardinalities and compression ratios, the cpu-times 
for the construction of the low-discrepancy sequence (cpu Halton seq.) and those for the computation of the compressed rules, where the new algorithm shows speed-ups from about 6 to more than 24 in the present degree range, ensuring moment residuals always below the required tolerance in at most 3 iterations. 
It is worth stressing a phenomenon already observed in \cite{ESV22}, that is possible failure of $Q_c^{dCATCH}$ which in some cases give much larger residuals than $Q_c^{bu}$.

\begin{figure}
    \centering
    \subfigure[]
        {\includegraphics[height=1.80in]{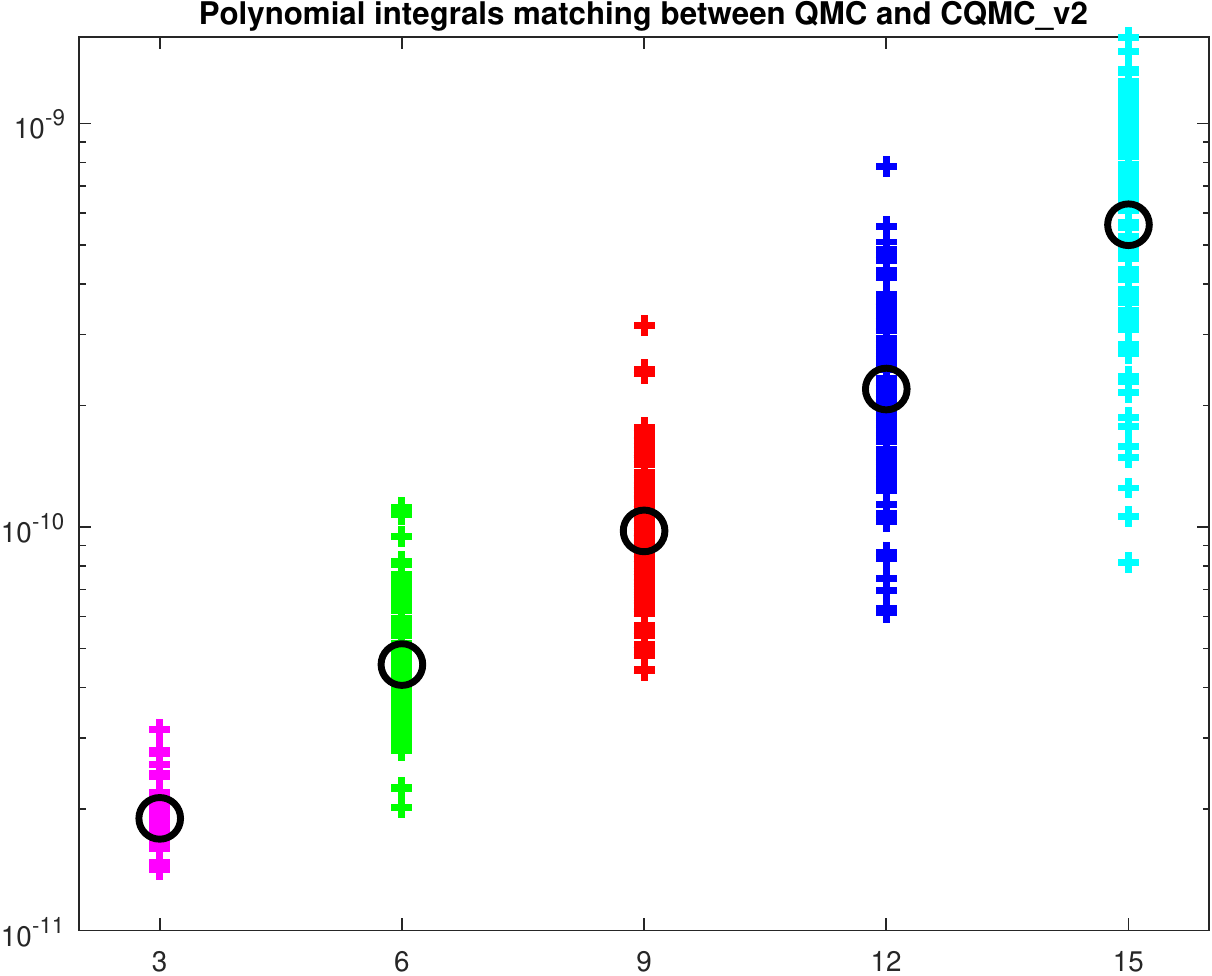}}
     ~ 
    \subfigure[]
        {\includegraphics[height=1.80in]{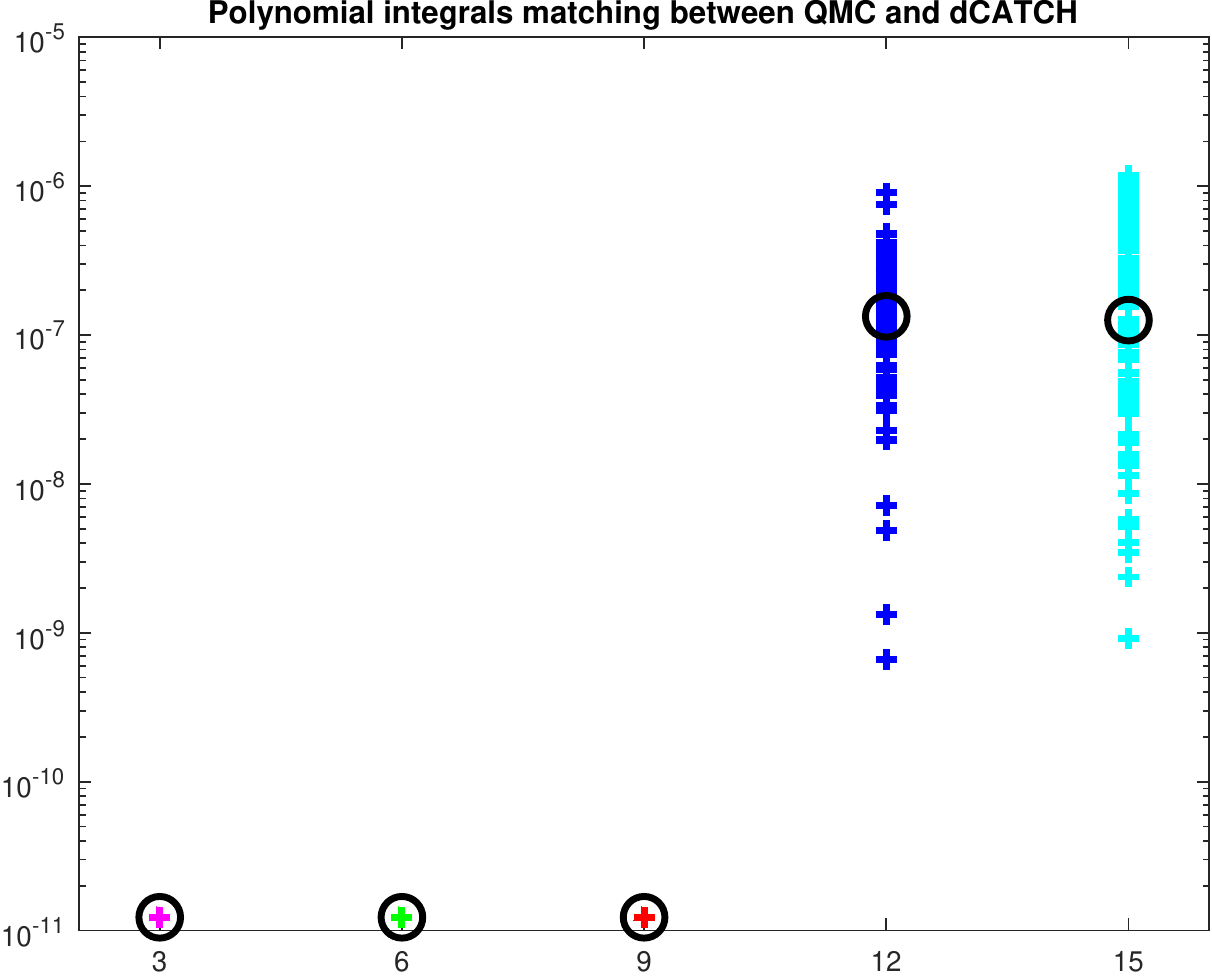}}
    \caption{  Relative QMC compression errors and their logarithmic average (circles) over 100 trials of random
polynomials for the bottom-up algorithm (left) and {\em{dCATCH}} (right) on the union of 3 balls. Note that the scales of the left and right figure are different. }
    \label{u3ball}
\end{figure}

\begin{figure}
    \centering   
    \subfigure[]
        {\includegraphics[height=1.80in]{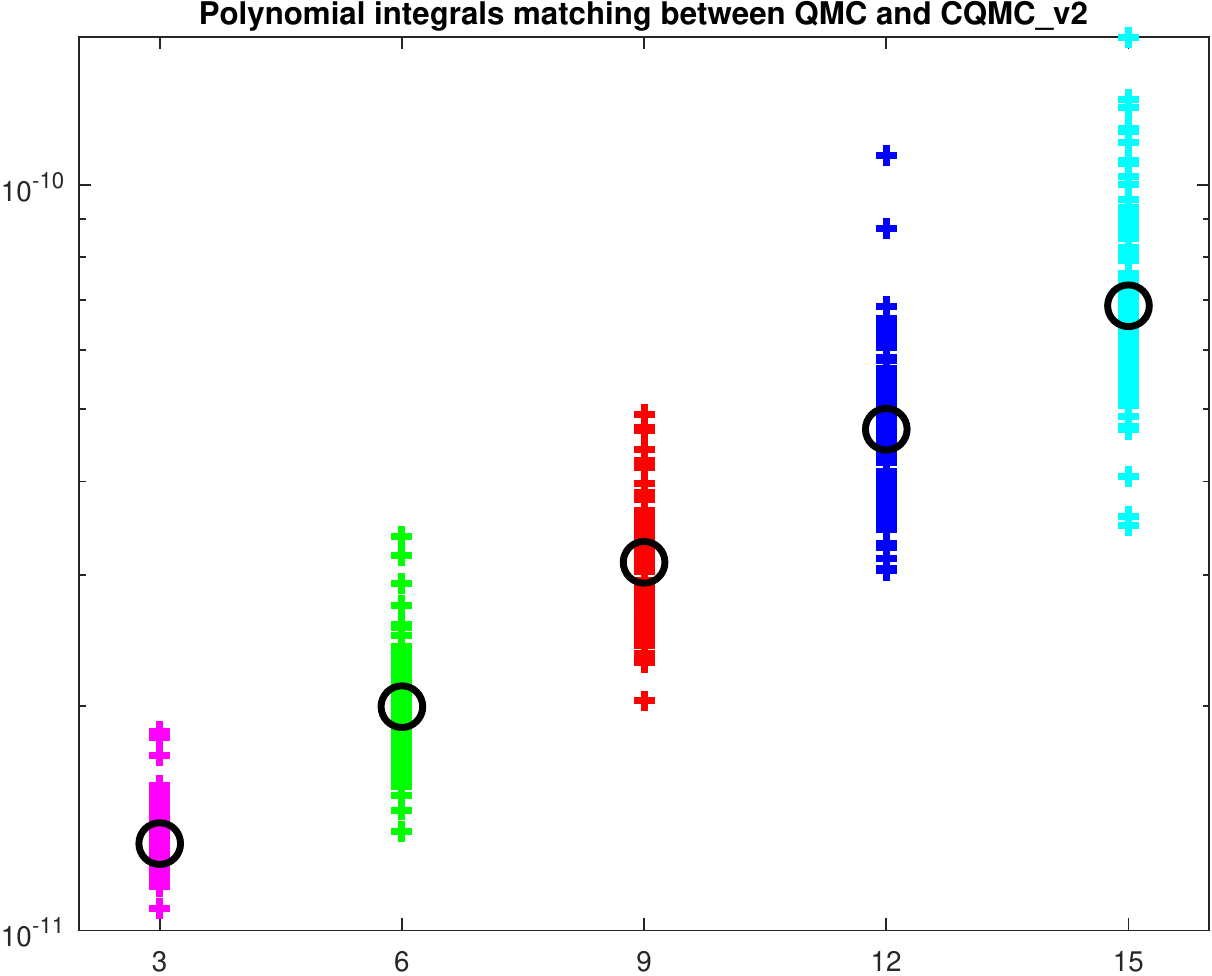}}
     ~ 
    \subfigure[]
        {\includegraphics[height=1.80in]{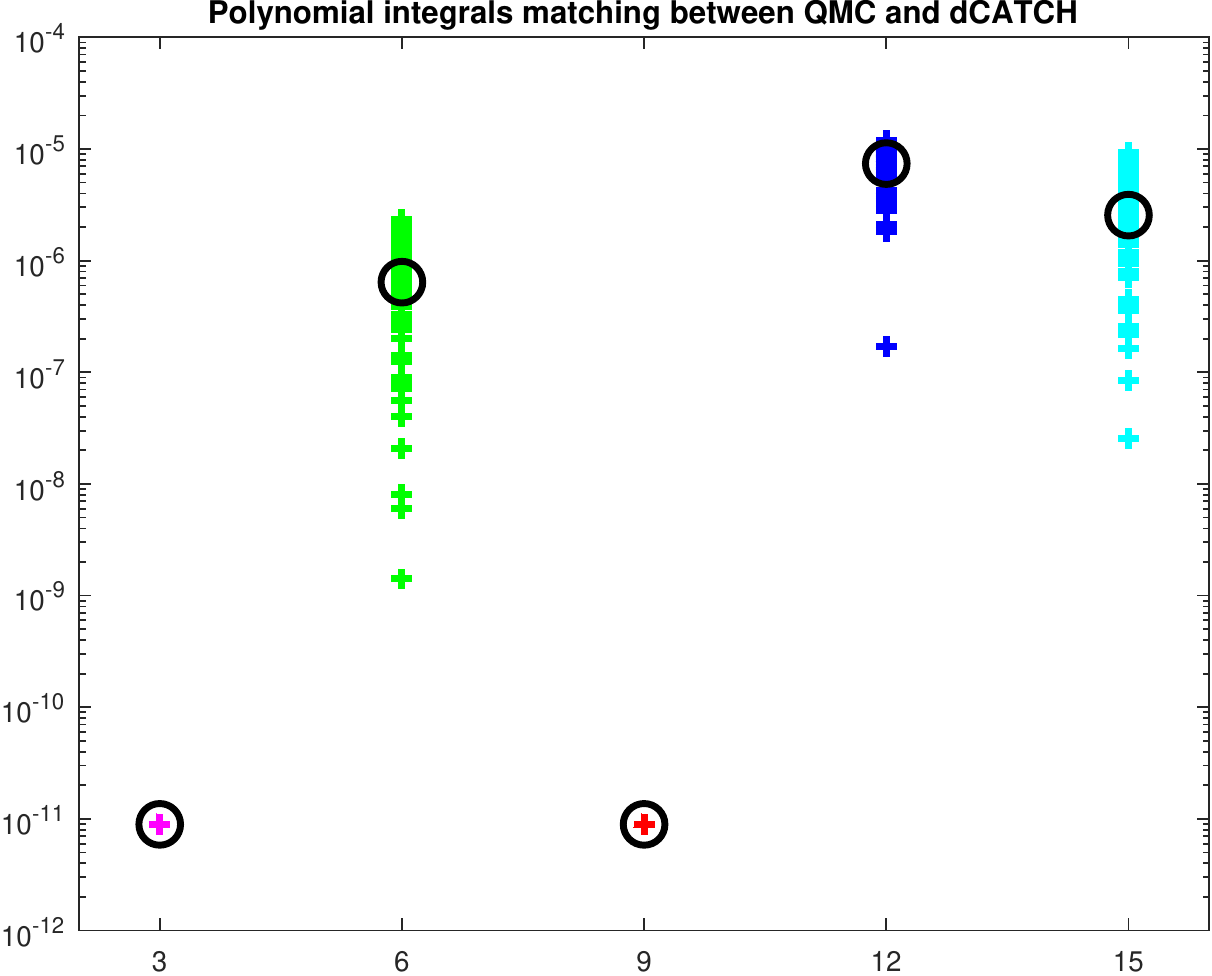}}
    \caption{  Relative QMC compression errors and their logarithmic average (circles) over 100 trials of random
polynomials for the bottom-up algorithm (left) and {\em{dCATCH}} (right) on the union of 100 balls. Note that the scales of the left and right figure are different. }
    \label{u100ball}
\end{figure}

In order to check polynomial exactness of the QMC compressed rules, in Figures {\ref{u3ball}}-{\ref{u100ball}} we show the relative QMC compression errors 
and their logarithmic averages (i.e. the sum of the log of the errors divided by the number of trials) over 100 trials of the polynomial 
\begin{equation} \label{polyex}
g(P)=(ax+by+cz+ d)^n\;,\;\;P=(x,y,z)
\end{equation}
where $a,b,c,d$ are uniform random variables in $[0,1]$. 
Moreover, in Tables {\ref{3balls-errors}}-{\ref{100balls-errors}} we show the integration {relative} errors on three test functions 
with different regularity, namely 
\begin{eqnarray} \label{testf}
f_1(P)&=&|P-P_0|^5 \nonumber\\
f_2(P)&=& \cos(x+y+z)\\
f_3(P)&=& \exp(-|P-P_0|^2) \nonumber
\end{eqnarray} 
where $P_0=(0,0,0)\in \Omega$, the first being of class $C^4$ with discontinuous fifth derivatives whereas the second and the third are  analytic. The reference values of the integrals have been computed by a QMC formula starting from $10^8$ Halton points in the bounding box. 

We see that the compressed formulas on more than one million points show errors 
of comparable order of magnitude, that as expected from estimate (\ref{equiv}) decrease while increasing the polynomial compression degree until they 
reach a size close to the QMC error   
(observe however that in Table {\ref{3balls-errors}} at degree $n=15$ only the bottom-up algorithm 
has reached the size of the QMC error).

\begin{table}[ht]
\begin{center}
{\footnotesize
\begin{tabular}{| c | c | c | c | c | c |}
\hline
deg & 3 & 6 & 9 & 12 & 15 \\
\hline \hline
card. $QMC$ & \multicolumn{5}{c|}{$M=$ 1,128,709}  \\
\hline
card. $Q_c^{dCATCH}$ & 20 & 84 & 220 & 452 & 806 \\
card. $Q_c^{bu}$    & 20 & 84 & 220 & 455 & 816 \\
compr. ratio & 5.6e+04 & 1.3e+04 & 5.1e+03 & 2.5e+03 & 1.4e+03 \\
\hline \hline  
cpu Halton seq. & \multicolumn{5}{c|}{9.0e-01s} \\ 
\hline 
cpu $Q_c^{dCATCH}$ & 3.4e+00s  & 1.9e+01s    & 4.9e+01s   & 1.4e+02s & 3.1e+02s  \\
\hline
cpu $Q_c^{bu}$ & 2.2e-01s  & 9.0e-01s   & 2.4e+00s  & 5.7e+00s & 2.6e+01s  \\
speed-up & 15.4 & 21.1 & 20.5 & 24.4 & 11.9 \\
\hline \hline
mom. resid. $Q_c^{dCATCH}$ & 8.9e-12 & 8.9e-12 & 8.9e-12 & $\star$ 5.1e-06 & $\star$ 1.1e-05\\
\hline
mom. resid. $Q_c^{bu}$ & & & & & \\
iter. 1 & 4.55e-16 & 1.51e-02 & 1.63e-01 & 3.81e-01 & 7.12e-01 \\
iter. 2 &         & 1.12e-15 & 1.85e-15 & 3.62e-15 & 8.06e-15\\
\hline
\end{tabular}
}
\caption{\small{Example with the union of 3 balls, in a bounding box with 2,400,000 low-discrepancy points.
}}
\label{3balls}
\end{center}
\end{table}

\begin{table}[ht]
\begin{center}
{\footnotesize
\begin{tabular}{| c | c | c | c | c | c |}
\hline
deg & 3 & 6 & 9 & 12 & 15 \\
\hline \hline
card. $QMC$ & \multicolumn{5}{c|}{$M=$ 1,195,806}  \\
\hline
card. $Q_c^{dCATCH}$ & 20 & 83 & 220 & 450 & 795 \\
card. $Q_c^{bu}$    & 20 & 84 & 220 & 455 & 816 \\
compr. ratio & 5.6e+04 & 1.3e+04 & 5.1e+03 & 2.8e+03  & 1.5e+03 \\
\hline \hline  
cpu Halton seq. & \multicolumn{5}{c|}{1.3e+00s}  \\ 
\hline 
cpu $Q_c^{dCATCH}$ & 3.4e+00s & 2.3e+01s  & 6.5e+01s    & 1.5e+02s & 3.7e+02s  \\
\hline
cpu $Q_c^{bu}$ & 2.5e-01s  & 8.7e-01s   & 2.6e+00s     & 9.5e+00s  & 6.7e+01s   \\
speed-up & 13.8 & 26.6 & 25.0 & 15.7 & 5.6 \\
\hline \hline
mom. resid. $Q_c^{dCATCH}$ & 1.1e-11 & $\star$ 1.2e-05 & 1.1e-11 & $\star$ 5.6e-05 & $\star$ 7.3e-05 \\
\hline
mom. resid. $Q_c^{bu}$ & & & & & \\
iter. 1 & 2.08e-16 & 9.41e-02 & 4.99e-01 & 1.51e+00 & 1.78e+00 \\
iter. 2 &   & 1.32e-15 & 2.20e-15 & 4.72e-15 & 8.30e-02 \\
iter. 3 &   &   &   &  & 7.32e-15 \\
\hline
\end{tabular}
}
\caption{\small{Example with the union of 100 balls, in a bounding box with 2,400,000 Halton points.
}}
\label{100balls}
\end{center}
\end{table}

\begin{table}[h]
\begin{center}
{\footnotesize
\begin{tabular}{| c | c | c | c | c | c |}
\hline
deg & 3 & 6 & 9 & 12 & 15 \\
\hline \hline
$E^{QMC}(f_1)$  & \multicolumn{5}{c|}{3.5e-04}  \\ 
\hline 
$E^{dCATCH}(f_1)$ &  1.3e-01  & 3.4e-04 & 3.5e-04 & 3.5e-04 & 3.5e-04 \\
$E^{bu}(f_1)$    & 2.3e-03 & 3.2e-04 & 3.5e-04 & 3.5e-04 & 3.5e-04 \\
 \hline \hline
$E^{QMC}(f_2)$  & \multicolumn{5}{c|}{7.3e-04}  \\ 
\hline  
$E^{dCATCH}(f_2)$ & 2.4e+00 & 7.0e-02 & 4.3e-03 & 7.3e-04 & 7.3e-04 \\
$E^{bu}(f_2)$    & 7.5e-01 & 3.7e-03 & 4.8e-04 & 7.4e-04 & 7.3e-04 \\
 \hline \hline
$E^{QMC}(f_3)$  & \multicolumn{5}{c|}{8.7e-05}  \\ 
\hline 
$E^{dCATCH}(f_3)$ & 7.1e-01 & 1.4e-01 & 9.4e-03 & 2.1e-03 & 1.1e-04\\
$E^{bu}(f_3)$    & 5.8e-01 & 2.8e-02 & 1.5e-02 & 9.5e-04 & 2.5e-05 \\
\hline 
\end{tabular}
}
\caption{\small{Example with 3 balls (the reference values are computed via QMC starting from $10^8$ Halton points in the bounding box).
}}
\label{3balls-errors}
\end{center}
\end{table}

\begin{table}[h]
\begin{center}
{\footnotesize
\begin{tabular}{| c | c | c | c | c | c |}
\hline
deg & 3 & 6 & 9 & 12 & 15 \\
\hline \hline
$E^{QMC}(f_1)$  & \multicolumn{5}{c|}{1.1e-04}  \\ 
\hline
$E^{dCATCH}(f_1)$ & 8.3e-02 & 8.8e-05 & 1.1e-04 & 1.1e-04 & 1.1e-04 \\
$E^{bu}(f_1)$    & 1.7e-03 & 9.8e-05 & 1.1e-04 & 1.1e-04 & 1.1e-04 \\
\hline \hline
$E^{QMC}(f_2)$  & \multicolumn{5}{c|}{1.7e-04}  \\ 
\hline 
$E^{dCATCH}(f_2)$ & 2.9e-01 & 8.7e-04 & 1.6e-04 & 1.7e-04 & 1.7e-04 \\
$E^{bu}(f_2)$    & 5.6e-02 & 1.5e-04 & 1.7e-04 & 1.7e-04 & 1.7e-04 \\
\hline \hline
$E^{QMC}(f_3)$  & \multicolumn{5}{c|}{2.2e-04}  \\ 
\hline
$E^{dCATCH}(f_3)$ & 2.3e-01 & 2.3e-03 & 8.4e-04 & 2.3e-04 & 2.2e-04\\
$E^{bu}(f_3)$    & 6.1e-03 & 3.6e-03 & 1.2e-04 & 2.3e-04 & 2.2e-04 \\
\hline 
\end{tabular}
}
\caption{\small{Example with 100 balls (the reference values are computed via QMC starting from $10^8$ Halton points in the bounding box).
}}
\label{100balls-errors}
\end{center}
\end{table}

\subsection{Surface integration on multibubbles}

We turn now to surface integration, on a domain $\Omega$ that is the boundary of an arbitrary union of balls, namely 
\begin{equation}\label{usph}
\Omega=\partial \bigcup_{j=1}^s B(C_j,r_j)=\bigcup_{j=1}^s\partial B(C_j,r_j)  {\backslash}  \bigcup_{j=1}^s \overset{\circ}{B} (C_j,r_j)\;,
\end{equation}
i.e. the set of all points lying on some sphere $\partial B(C_j,r_j)$, $j=1,\ldots,s$, but not internally to any of the balls $ B(C_{k},r_{k})$, ${k} \neq j$. We present two examples, corresponding to the same centers and radii considered above for volume integration, i.e. the surface of the union of 3 balls and of 100 balls in Section 3.1. Notice that $\Omega$ is a subset of an algebraic surface, i.e. the union of the corresponding spheres. Though the polynomial spaces dimension could be computed theoretically by algebraic 
geometry methods (cf., e.g., \cite{CLOS15}), we do not enter this delicate matter here, since the algorithm computes numerically such a dimension by a rank revealing approach on a Vandermonde-like matrix. 

In this case we have applied the extension discussed in Remark {\ref{QMC-union-seq}}, constructing an evenly distributed sequence $X_M$ on the whole $\Omega$ by taking a large number of low discrepancy points on each sphere $\partial B(C_j,r_j)$, and then selecting those belonging to the portions of the sphere that contribute to the surface of the union, that are those not internal to any other ball. Namely, we have taken on each sphere the mapped Halton points from the rectangle $[-1,1] \times [0,2\pi]$ by the {\em area preserving} transformation 
\begin{equation} \label{mapped}
(t,\phi) \mapsto C_j+r_j(\sqrt{1-t^2}\,\cos(\phi),\sqrt{1-t^2}\,\sin(\phi),t)\;,
\end{equation}
which preserves also the 
low-discrepancy property. The points are finally ordered by picking alternatively one point per active portion 
of the surface of the union, with a local weight attached to each point. An illustration of compressed points extracted starting from 4000 mapped Halton points on each sphere is given in Figure \ref{compressed-deg9}.

In Tables {\ref{3spheres}}-{\ref{100spheres}} we report for this surface integration examples the same quantities appearing in Tables {\ref{3balls}}-{\ref{100balls} for the volume integration, where we use again the $dCATCH$ code in \cite{DMV20} to compress the QMC formula on the whole $X_M$, since also that algorithm was conceived to work with polynomial spaces possibly restricted to algebraic surfaces. Here we start from 500,000 mapped Halton points on each sphere in the 3 balls example, and from 60,000 in the 100 balls instance, obtaining a sequence of about one million low-discrepancy points on the corresponding ball union surfaces. As before we set $m_{j+1} = 2m_j$, $j\geq 1$ with $m_1 = 2N$ and $\varepsilon=10^{-10}$.

Again we get impressive compression ratios, and speed-ups varying from about 5 to more than 16. Moreover, the bottom-up algorithm gives always a residual 
below the given tolerance, whereas $dCATCH$ turns out to be more prone to failure (see the 
residuals for degree $n=15$ in the example with 3 balls and degrees $n=9,15$ in the example with 100 balls). 

The logarithmic average errors concerning surface integration of the random polynomial (\ref{polyex}), restricted to the boundary 
of the union, are plotted in Figures {\ref{u3sph}}-{\ref{u100sph}}. In Tables {\ref{3spheres-errors}}-{\ref{100spheres-errors}} we show the surface integration errors for the three test functions in (\ref{testf}), where $P_0$ is a suitably chosen point on the surface of the ball union. We see again that the compressed formulas on more than one million points show errors 
of comparable order of magnitude, that as expected from estimate (\ref{equiv}) decrease while increasing the polynomial compression degree, until they 
reach a size close to the QMC error.   

\begin{table}[ht]
\begin{center}
{\footnotesize
\begin{tabular}{| c | c | c | c | c | c |}
\hline
deg & 3 & 6 & 9 & 12 & 15\\
\hline \hline
card. $QMC$ & \multicolumn{5}{c|}{$M=$ 1,024,179}  \\
\hline
card. $Q_c^{dCATCH}$ & 20 & 83 & 200 & 371 & 572 \\
card. $Q_c^{bu}$    & 20 & 83 & 200 & 371 & 596 \\
compr. ratio & 5.1e+04 & 1.2e+04 & 5.1e+03 & 2.8e+03 & 1.7e+03  \\
\hline \hline  
cpu Halton seq. & \multicolumn{5}{c|}{8.8e-01s}  \\ 
\hline 
cpu $Q_c^{dCATCH}$ & 2.8e+00s & 1.7e+01s & 5.0e+01s & 1.4e+02s  & 3.2e+02s  \\
\hline
cpu $Q_c^{bu}$ & 3.1e-01s & 1.1e+00s  & 2.7e+00s   & 5.8e+00s   & 6.5e+01s   \\
speed-up & 9.0 & 15.9 & 18.2 & 24.0 & 4.9 \\
\hline \hline
mom. resid. $Q_c^{dCATCH}$ & 8.6e-12 & 8.9e-12 & 8.9e-12 & 8.9e-12 & $\star$ 1.6e-06 \\
\hline
mom. resid. $Q_c^{bu}$ & & & & & \\
iter. 1 & 7.2e-01 & 1.4e-15 & 2.8e-15 & 4.2e-15 & 2.6e-01 \\
iter. 2 & 3.7e-16 & & & & 1.3e-01\\
iter. 3 & & & & &  3.3e-12  \\
\hline
\end{tabular}
}
\caption{\small{Compression of surface QMC integration on the union of 3 balls, starting from 500,000 low-discrepancy points on each sphere.
}}
\label{3spheres}
\end{center}
\end{table}

\begin{table}[ht]
\begin{center}
{\footnotesize
\begin{tabular}{| c | c | c | c | c | c |}
\hline
deg & 3 & 6 & 9 & 12 & 15 \\
\hline \hline
card. $QMC$ & \multicolumn{5}{c|}{$M=$ 1,032,718}  \\
\hline
card. $Q_c^{dCATCH}$ & 20 & 84 & 219 & 455 & 807 \\
card. $Q_c^{bu}$    & 20 & 84 & 220 & 455 & 816 \\
compr. ratio & 5.2e+04 & 1.2e+04 & 4.7e+03 & 2.3e+03 & 1.3e+03 \\
\hline \hline  
cpu Halton seq. & \multicolumn{5}{c|}{1.5e+01s}  \\ 
\hline 
cpu $Q_c^{dCATCH}$ & 2.8e+00s & 1.6e+01s & 4.3e+01s   & 1.1e+02s & 2.4e+02s  \\
\hline
cpu $Q_c^{bu}$ & 3.2e-01s   & 1.1e+00s    & 3.0e+00s   & 6.8e+00s & 2.4e+01s \\
speed-up & 8.7 & 14.5 & 14.3 & 16.2 & 9.9 \\
\hline \hline 
mom. resid. $Q_c^{dCATCH}$ & 9.0e-13 & 9.1e-13 & $\star$ 3.2e-06 & 9.3e-13 & $\star$ 1.8e-05 \\
\hline
mom. resid. $Q_c^{bu}$ & & & & & \\
iter. 1 & 2.09e+00 & 1.22e+00 & 5.53e-01 & 6.18e-01 & 1.41e-01 \\
iter. 2 & 7.49e-16 & 1.30e-15 & 2.52e-15 & 5.16e-15 & 1.18e-14 \\
\hline
\end{tabular}
}
\caption{\small{Compression of surface QMC integration on the union of 100 balls, starting from 60,000 low-discrepancy points on each sphere.
}}
\label{100spheres}
\end{center}
\end{table}

\begin{figure}
    \centering
    \subfigure[]
        {\includegraphics[height=1.80in]{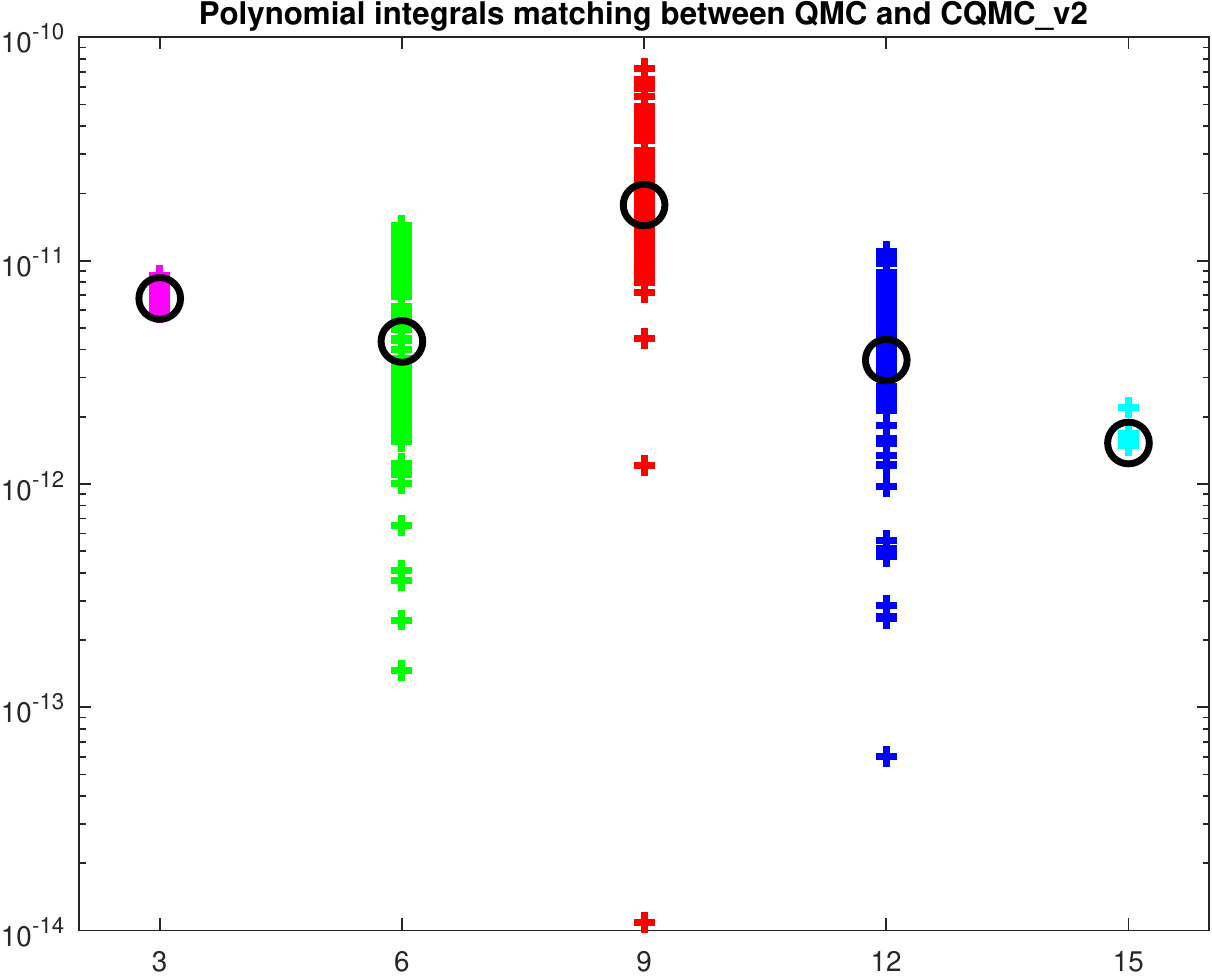}}
     ~ 
    \subfigure[]
        {\includegraphics[height=1.80in]{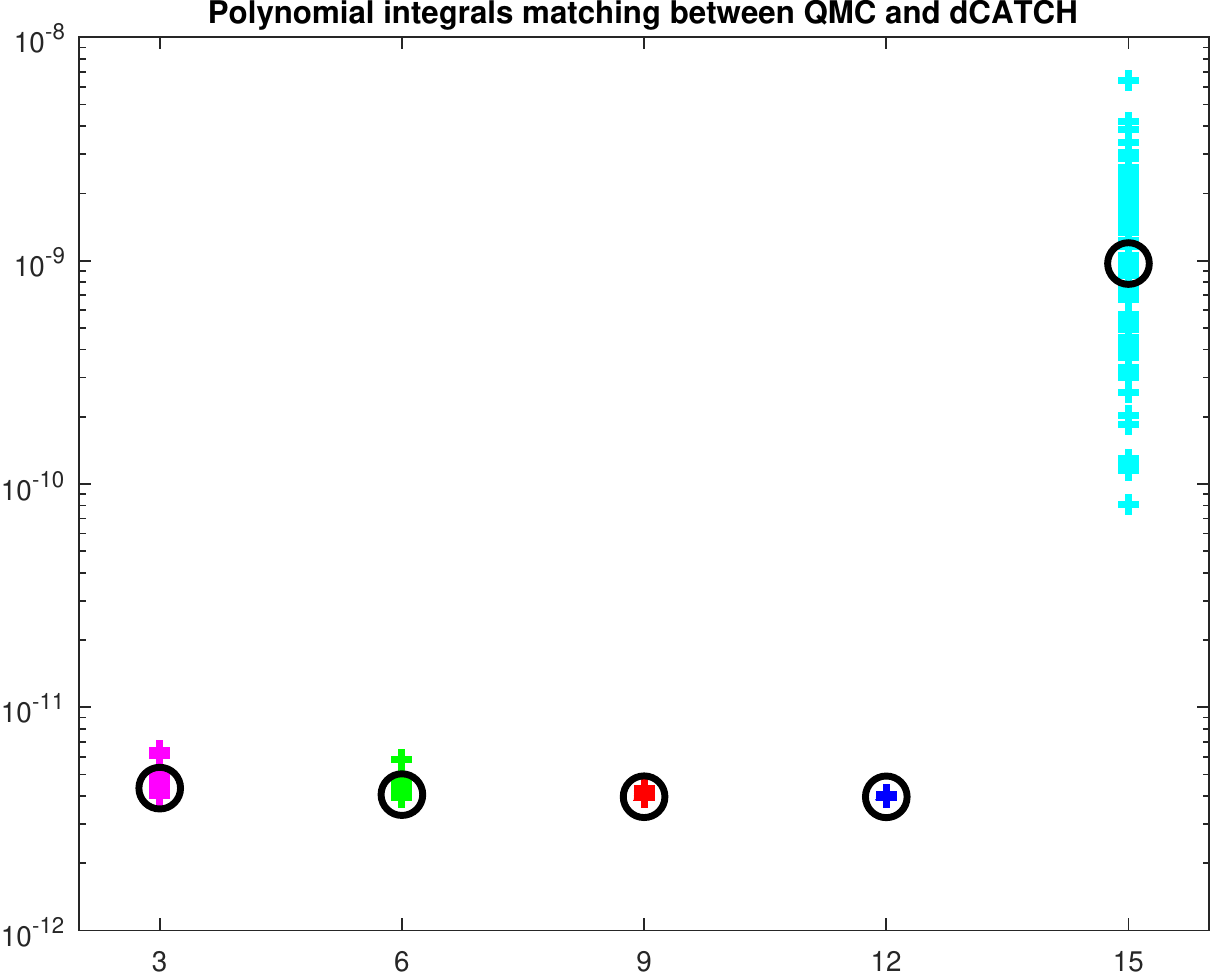}}
    \caption{ Relative QMC compression errors and their logarithmic average (circles) over 100 trials for the bottom-up algorithm (left) and {\em{dCATCH}} (right) of random
polynomials on the surface of the union of 3 balls.  Note that the scales of the left and right figure are different. }
    \label{u3sph}
\end{figure}

\begin{figure}
    \centering
    \subfigure[]
        {\includegraphics[height=1.80in]{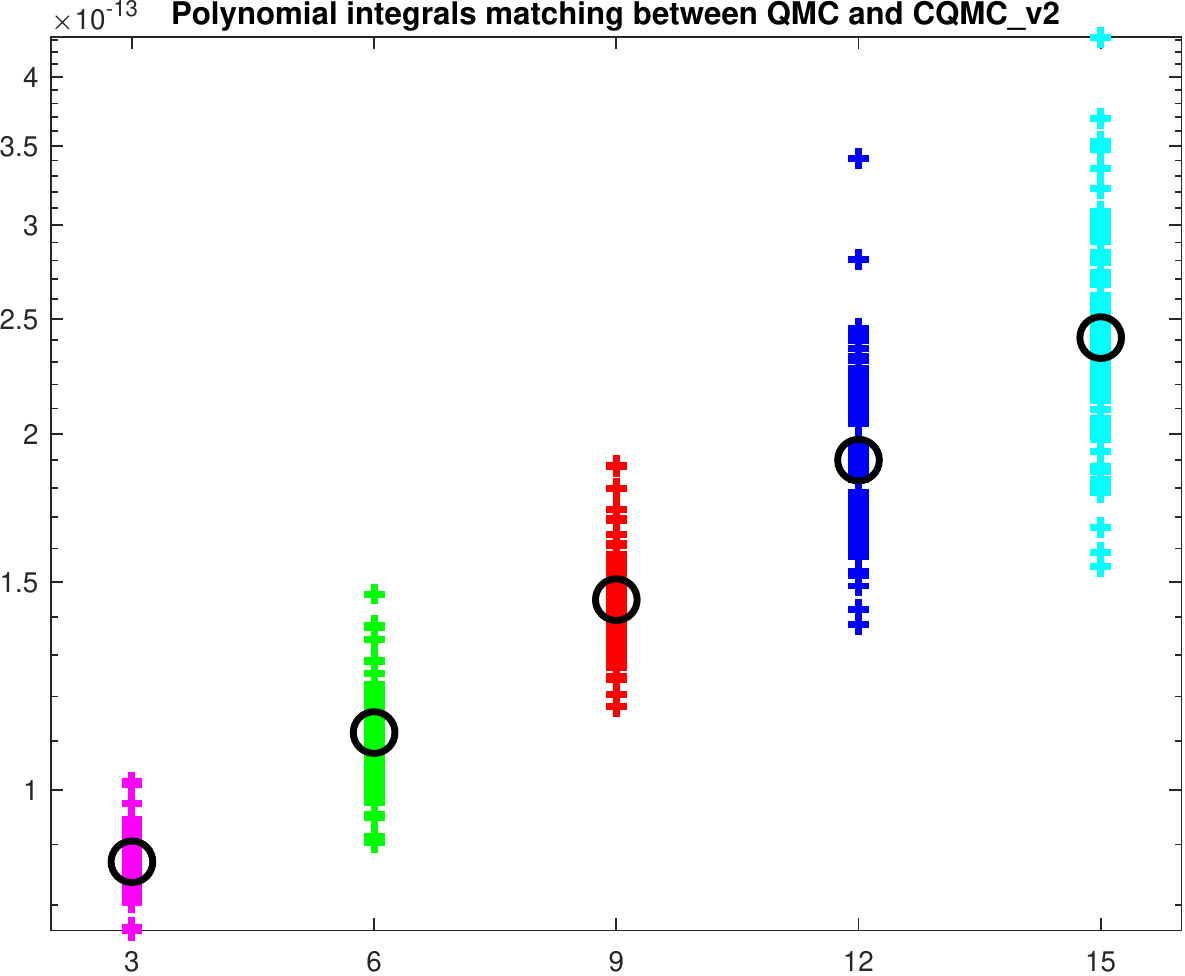}}
     ~ 
    \subfigure[]
        {\includegraphics[height=1.80in]{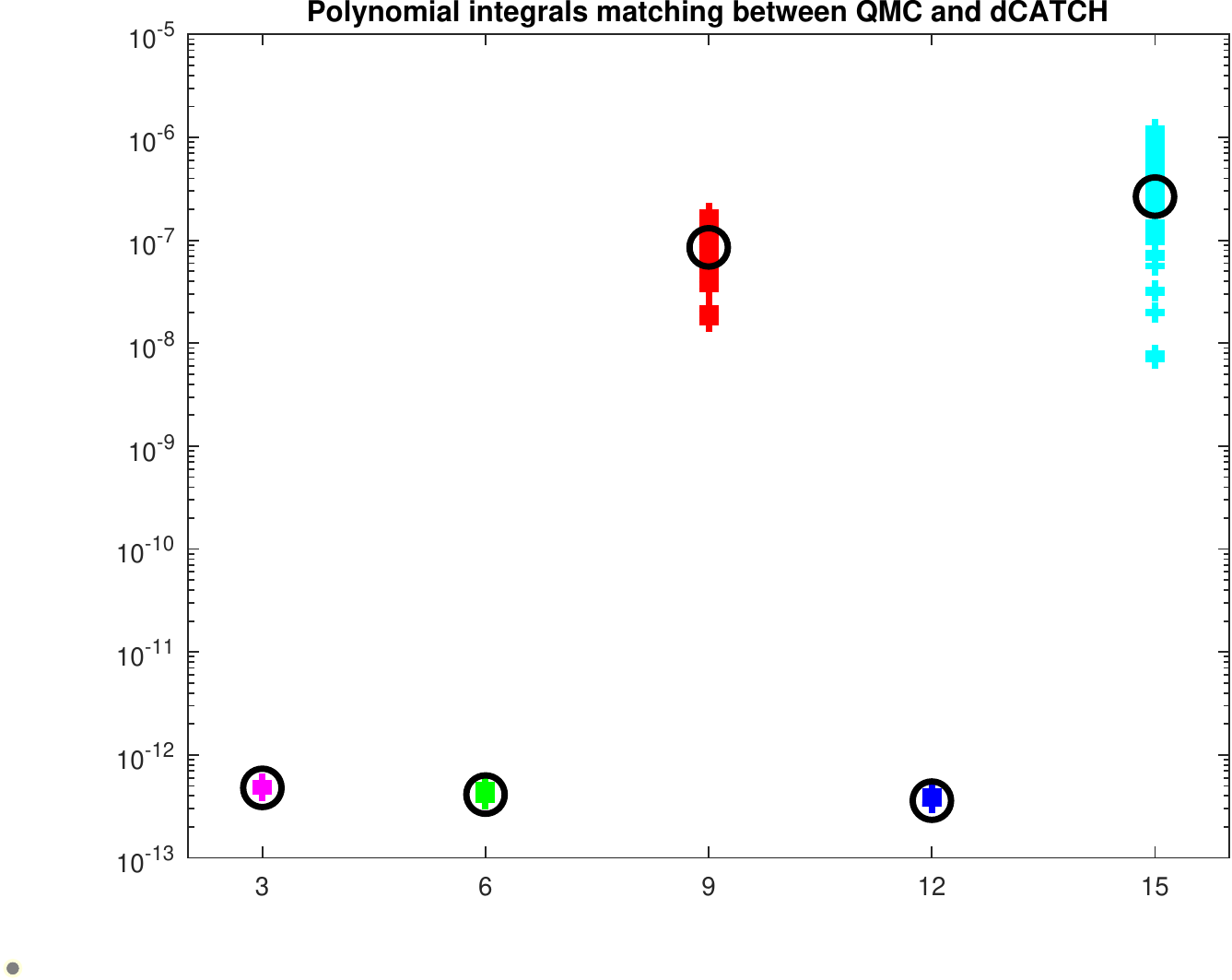}}
    \caption{ Relative QMC compression errors and their logarithmic average (circles) over 100 trials for the bottom-up algorithm (left) and {\em{dCATCH}} (right) of random polynomials on the surface of the union of 100 balls. Note that the scales of the left and right figure are different.}
    \label{u100sph}
\end{figure}

 \begin{table}[h]
\begin{center}
{\footnotesize
\begin{tabular}{| c | c | c | c | c | c |}
\hline
deg & 3 & 6 & 9 & 12 & 15 \\
\hline \hline
$E^{QMC}(f_1)$  & \multicolumn{5}{c|}{3.9e-06}  \\ 
\hline 
$E^{dCATCH}(f_1)$ & 2.3e-04 & 8.3e-06 & 4.0e-06 & 3.9e-06 & 3.9e-06 \\
$E^{bu}(f_1)$    & 3.7e-04 & 3.6e-06 & 4.0e-06 & 3.9e-06 & 3.9e-06\\
\hline 
$E^{QMC}(f_2)$  & \multicolumn{5}{c|}{8.6e-05}  \\ 
\hline 
$E^{dCATCH}(f_2)$ & 3.5e-01 & 3.2e-02 & 8.3e-04 & 8.4e-05 & 8.6e-05\\
$E^{bu}(f_2)$    & 1.3e+00 & 2.2e-02 & 8.3e-06 & 8.5e-05 & 8.6e-05 \\
\hline 
$E^{QMC}(f_3)$  & \multicolumn{5}{c|}{5.8e-06}  \\ 
\hline 
$E^{dCATCH}(f_3)$ & 3.9e-01 & 5.8e-03 & 6.9e-04 & 5.8e-05 & 8.9e-06\\
$E^{bu}(f_3)$    & 2.5e-02 & 3.8e-03 & 4.7e-06 & 9.1e-05 & 6.0e-06\\
\hline 
\end{tabular}
}
\caption{\small{Compression of surface QMC integration on the union 3 balls (the reference values are computed via QMC starting from $10^6$ points on each sphere).
}}
\label{3spheres-errors}
\end{center}
\end{table}

\begin{table}[h]
\begin{center}
{\footnotesize
\begin{tabular}{| c | c | c | c | c | c |}
\hline
deg & 3 & 6 & 9 & 12 & 15 \\
\hline \hline
$E^{QMC}(f_1)$  & \multicolumn{5}{c|}{4.0e-05}  \\ 
\hline 
$E^{dCATCH}(f_1)$ & 2.9e-03 & 2.8e-05 & 3.9e-05 & 4.0e-05 & 4.0e-05 \\
$E^{bu}(f_1)$    & 3.2e-02 & 3.7e-05 & 3.9e-05 & 4.0e-05 & 4.0e-05 \\
\hline 
$E^{QMC}(f_2)$  & \multicolumn{5}{c|}{2.0e-04}  \\ 
\hline 
$E^{dCATCH}(f_2)$ & 1.4e-01 & 6.4e-04 & 1.7e-04 & 2.0e-04 & 2.0e-04 \\
$E^{bu}(f_2)$    & 1.4e-01 & 6.2e-05 & 1.9e-04 & 2.0e-04 & 2.0e-04 \\
\hline 
$E^{QMC}(f_3)$  & \multicolumn{5}{c|}{1.6e-04}  \\ 
\hline 
$E^{dCATCH}(f_3)$ & 2.9e-02 & 1.7e-02 & 4.3e-04 & 1.5e-04 & 1.6e-04\\
$E^{bu}(f_3)$    & 1.6e-02 & 1.7e-03 & 2.6e-04 & 1.5e-04 & 1.6e-04 \\
\hline 
\end{tabular}
}
\caption{\small{Compression of surface QMC integration on the union 100 balls (the reference values are computed via QMC starting from $10^6$ points on each sphere).
}}
\label{100spheres-errors}
\end{center}
\end{table}

\section{Acknowledgements}
Work partially
supported by the
DOR funds and the biennial project BIRD 192932
of the University of Padova, and by the INdAM-GNCS 
2022 Project ``Methods and software for multivariate integral models''.
This research has been accomplished within the RITA ``Research ITalian network on Approximation", the UMI Group TAA ``Approximation Theory and Applications" (G. Elefante, A. Sommariva) and the SIMAI Activity Group ANA\&A (A. Sommariva, M. Vianello).

\end{document}